\documentstyle[epsf,amsfonts,addeng]{article}
\author{Yu.~M.~Burman\thanks{Independent University of Moscow; e-mail:
burman@mccme.ru}}
\title{Relative moduli spaces of complex structures: an example}
\date{}

\def \Diff {\name{Diff}}
\def \Map {\name{Map}}
\def \Mod {{\cal M}}
\def \Equiv {\name{\cal E}}
\def \Proj {\name{\cal P}}
\def \Lift {\name{\cal L}}
\def \Compl {\name{\cal J}}
\def \Teich {{\cal T}}
\def \Group {{\frak g}}

\def \a {\alpha}
\def \b {\beta}
\def \one {^{(1)}}
\def \two {^{(2)}}

\begin{document}
\maketitle
\begin{abstract}
Let  $M$ and $N$ be even-dimensional oriented real manifolds, and
$u:M \to N$ be a smooth mapping. A pair of complex structures at $M$ and
$N$ is called $u$-compatible if the mapping $u$ is holomorphic with respect
to these structures. The quotient of the space of $u$-compatible pairs of
complex structures by the group of $u$-equivariant pairs of diffeomorphisms
of $M$ and $N$ is called a moduli space of $u$-equivariant complex
structures. The paper contains a description of the fundamental group $G$
of this moduli space in the following case: $N = CP^1$, $M \subset CP^2$ is
a hyperelliptic genus $g$ curve given by the equation $y^2 = Q(x)$ where
$Q$ is a generic polynomial of degree $2g+1$, and $u(x,y) = y^2$. The group
$G$ is a kernel of several (equivalent) actions of the braid-cyclic group
$BC_{2g}$ on $2g$ strands. These are: an action on the set of trees with
$2g$ numbered edges, an action on the set of all splittings of a
$(4g+2)$-gon into numbered nonintersecting quadrangles, and an action on a
certain set of subgroups of the free group with $2g$ generators. $G_{2g}
\subset BC_{2g}$ is a subgroup of the index $(2g+1)^{2g-2}$.

{\em Key words}: Teichm\"uller spaces, Lyashko--Looijenga map, braid group.
\end{abstract}

\section*{Introduction}

Let two even-dimensional oriented real manifolds $M^{2m}$ and $N^{2n}$ be
given, and $u:M \to N$ be a smooth mapping. A pair of complex structures at
$M$ and $N$ will be called {\em $u$-compatible} if the mapping $u$ is
holomorphic with respect to these structures. The set of all $u$-compatible
pairs of complex structures will be denoted $\Compl(u)$. Call a pair of
diffeomorphisms $\a: M \to M$ and $\b: N \to N$ {\em $u$-equivariant} if
the diagram
\begin{equation}\label{Eq:Equivar}
\begin{array}{lcl}
M & {\stackrel{\mbox{\scriptsize $\a$}}{\longrightarrow}} & M \\
\downarrow u && \downarrow u\\
N & {\stackrel{\mbox{\scriptsize $\b$}}{\longrightarrow}} & N
\end{array}
\end{equation}
commutes. The $u$-equivariant pairs of diffeomorphisms form a
topological group $\Equiv(u)$. The connected component of the unity
in the group $\Equiv(u)$ consists of pairs of diffeomorphisms
homotopic to identity (inside $\Equiv(u)$). This subgroup will be
denoted $\Equiv_0(u)$.

Groups $\Equiv(u)$ and $\Equiv_0(u)$ act on the space $\Compl(u)$. A
quotient $\Teich(u) = \Compl(u) / \Equiv_0(u)$ will be called a
Teichm\"uller space of $u$-compatible pairs of complex structures. The
smaller quotient $\Mod(u) = \Compl(u) / \Equiv(u)$ will be called a moduli
space. These spaces are ``relative'' analogs of the usual Teichm\"uller
space and moduli space of complex structures reducing to them in the case
when $u: M \to N$ is a diffeomorphism.

The article contains a description of the fundamental group
$\pi_1(\Mod(u))$ in the following case. The target manifold $N$ is
$\Complex P^1$.  The manifold $M \subset \Complex P^2$ is a sphere with $g$
handles (denoted $M_g$) given in homogeneous coordinates as $\{[x:y:z] \mid
y^2z^{2g-1} = q(x,z)\}$; here $q(x,z)$ is a homogeneous polynomial of
degree $2g+1$ such that all its critical values are different. The mapping
$u$ is defined by the formula $u([x:y:z]) = [y^2:z^2]$.

The article has the following structure. In Section~\ref{Sec:Fund} one
shows (Theorem~\ref{Th:Pi1IsMap}) that $\pi_1(\Mod(u)) = \Map(u) \bydef
\Equiv(u) / \Equiv_0(u)$ provided $M$ and $N$ have real dimension $2$, $u$
is a covering with at least $3$ branching values, and there exist complex
structures on $M$ and $N$ such that $u$ is holomorphic (i.e., $\Mod(u) \ne
\emptyset$). This result extends the corresponding theorem of the
``absolute'' Teichm\"uller theory and reduces the problem to the study of
the relative mapping class group $\Map(u)$. In Section~\ref{Sec:Equiv} we
prove (Theorem~\ref{Th:ExaSeq}) that if $u:M_g \to \Complex P^1$ is the
mapping described above then the group $\Map(u)$ is isomorphic, up to a
$\Integer_2$-extension, to its projection $\Lift(u)$ onto the group
$G_{2g}$ of diffeomorphisms of $\Complex P^1$ preserving branching values
of $u$. In Subsection~\ref{SSec:Lift} it is also proved
(Theorem~\ref{Th:LiftStab}) that $\Lift(u)$ is isomorphic to a point
stabilizer of a certain action of the braid-cyclic group $BC_{2g}$ on the
set of all splittings of a $(4g+2)$-gon into nonintersecting quadrangles.

In Section~\ref{Sec:Action} this action is studied in detail. An algorithm
for computing the action is given in Subsection~\ref{SSec:ActQuad}. In
Subsection~\ref{SSec:ActTree} we define an action of $BC_N$ on another
object, a set of trees with $N$ edges marked $0,1,\dots,N-1$, and show
that this two actions are related by some geometric construction. It is
proved that the action of $BC_N$ on trees is transitive, which allows to
obtain some information about the group $\Lift(u)$ (Corollary of
Statement~\ref{St:TreeToQuad}). We formulate a conjecture about
presentation of the group $\Lift(u)$. Our construction has some relations
with the theory of Lyashko--Looijenga mapping; these are also considered in
Subsection~\ref{SSec:ActTree}. The last Subsection~\ref{SSec:Subgr} is a
sort of appendix, it is devoted to the algebraic ``roots'' of the action
considered above.

\subsection*{Acknowledgments}

The initial inspiration of this work was an article \cite{AniLan} of
S.~Anisov and S.~Lando considering in fact the genus $1$ case of the
problem (though Teichm\"uller spaces do not appear directly in the
article). Several important errors in the earlier versions of the text were
pointed out by V.~Arnold, M.~Kontsevich, A.~Zvonkin, and D.~Zvonkine. The
author thanks B.~Wajnryb for making him acquainted with the paper
\cite{Thurston}. The author is grateful to L.~Funar for his letter
concerning difficulties in construction of the homomorphisms from the braid
groups to $\Map(M_g)$. An elegant proof of the lemma at
page~\pageref{Mk:Autom} was suggested by V.~Prasolov. Author is also
grateful to F.~Cohen for fruitful discussions.

Most of the work was completed during the author's stay in the University
Paris IX supported by a grant of the French Ministry of Foreign Affairs ---
the author is using an opportunity to express his warm gratitude to both
institutions.

\section{Fundamental group of the moduli space}\label{Sec:Fund}

\subsection{The main players}\label{SSec:Def}

For convenience we summarize here the definitions of the main objects of
the article. Some of them were already described in the Introduction. The
reader should refer to the list when necessary.

\begin{itemize}
\item $M$ and $N$ --- orientable $2$-dimensional (real) manifolds.
\item $u:M \to N$ --- a ramified covering, i.e.\ a smooth mapping such that
the preimage $u^{-1}(x) \subset M$ of any point $x \in N$ is discrete. The
set of ramification points of $u$ is denoted $R(u) \subset M$, and the set
of branching values (images of ramification points), $B(u) \subset N$.
\item $\Equiv(u)$ --- a topological group of $u$-equivariant pairs of
diffeomorphisms, i.e.\ of the pairs of diffeomorphisms $\a:M \to M$, $\b:N
\to N$ such that diagram~(\ref{Eq:Equivar}) commutes.
\item $\Equiv_0(u) \subset \Equiv(u)$ --- a connected component of the
unity.
\item $\Map(u) = \Equiv(u)/\Equiv_0(u)$ --- the $u$-equivariant mapping
class group.
\item $\Compl(u)$ --- the set of $u$-compatible pairs of complex
structures, i.e.\ pairs of complex structures at $M$ and $N$ such that the
mapping $u$ is holomorphic with respect to them.
\item $\Teich(u) = \Compl(u)/\Equiv_0(u)$, $\Mod(u) = \Compl(u)/\Equiv(u) =
\Teich(u)/\Map(u)$ --- relative Teichm\"uller space and moduli space.
\item $\Diff(K)$ --- the set of diffeomorphisms of the manifold $K$ (most
often $K = M$ or $N$); $\Diff(K;S_1, \dots, S_n)$ --- the set of
diffeomorphisms of the manifold $K$ preserving subsets $S_1, \dots, S_n
\subset K$ (each one individually, but not necessarily pointwise);
$\Diff_0(K)$, $\Diff_0(K;S_1, \dots, S_n)$ --- corresponding connected
components of the unity; $\Map(K) = \Diff(K)/\Diff_0(K)$, $\Map(K;S_1,
\dots, S_n) = \Diff(K;S_1, \dots, S_n)/\Diff_0(K;S_1, \dots, S_n)$ ---
corresponding mapping class groups.
\item $p_1, p_2$ --- generic names for mappings relating to a pair its
first, resp., second, term. These mappings will be applied to the pairs of
complex structures, diffeomorphisms, mapping classes, etc.
\item $\Equiv\one(u) \bydef p_1(\Equiv(u)) \subset \Diff(M,R(u)) \subset
\Diff(M)$ (apparently a $u$-equivariant mapping sends ramification points
to ramification points), $\Equiv\two(u) \bydef p_2(\Equiv(u)) \subset
\Diff(N,B(u)) \subset \Diff(N)$ (the same is true for branching values).
$\Equiv_0\one(u) \bydef p_1(\Equiv_0(u))$, $\Equiv_0\two(u) \bydef
p_2(\Equiv_0(u))$.
\item The quotient $\Equiv\one(u)/\Equiv_0\one(u) \subset \Map(M,R(u))$
will be called a group of projectable mapping classes (at $M$); it will be
denoted $\Proj(u)$. Similarly, $\Lift(u) \bydef
\Equiv\two(u)/\Equiv_0\two(u) \subset \Map(N,B(u))$ will be called a group
of liftable mapping classes (at $N$).
\end{itemize}

\subsection{Homotopy structure}\label{SSec:Homot}

Our aim in this Subsection is the following theorem:

\begin{theorem}\label{Th:Pi1IsMap}
Let $B(u)$ contain at least three points. Then $\pi_1(\Mod(u))$ is
isomorphic to the relative mapping class group $\Map(u) \bydef \Equiv(u) /
\Equiv_0(u)$.
\end{theorem}

To prove it, we are to prove first some auxiliary statements about homotopy
structure of diffeomorphism groups and Teichm\"uller spaces. Almost
everything we do is based on the classical lemma due to Alexander:

\begin{Lemma}
The group of orientation-preserving diffeomorphisms of the $n$-dimensional
disk with pointwise fixed boundary is contractible.
\end{Lemma}
See \cite{Alexander} for proof.

\begin{Corollary}
Let $S$ be a noncompact $2$-manifold, and $a, b \in S$. Then any connected
component of the space of paths in $S$ beginning at $a$ and ending in $b$
is contractible.
\end{Corollary}

\begin{proof}
The manifold $S$ can be contracted to its $1$-dimensional subcomplex
(graph) such that $a$ and $b$ are its vertices. A connected component in
the space of paths joining $a$ and $b$ in this graph is obviously
contractible.
\end{proof}

\begin{statement}\label{St:E0Contr}
Let the set $B(u)$ of branching values contain at least three points. Then
the set $\Equiv_0(u)$ is contractible.
\end{statement}

\begin{proof}
Since $u$ is a ramified covering, the mapping $p_2: \Equiv_0(u) \to
\Equiv_0\two(u)$ is a covering, too (nonramified). So it suffices to prove
that the set $\Equiv_0\two(u)$ is contractible.

The set $B(u) \subset N$ is discrete. Join its points with a {\em network},
a set of smooth nonintersecting arcs such that the complement of these arcs
in $M$ is a union of discs. Apparently a liftable diffeomorphism homotopic
to the identity maps each branching value to itself. A network is mapped to
another network homotopic to the original one. Thus $\Equiv_0\two(u)$ is
fibered over the connected component of the space of networks. The fiber is
a direct product of several copies of the group of diffeomorphisms of a
$2$-disk, fixed on the boundary. The fiber is contractible by Alexander's
lemma, and the base is contractible by its Corollary given in the beginning
of this Subsection.
\end{proof}

\begin{Corollary}
Spaces $\Equiv_0(u)$ and $\Equiv_0\two(u)$ are homeomorphic (the covering
$p_2$ is trivial).
\end{Corollary}

\begin{statement}\label{St:TContr}
Let the set $B(u)$ contain at least three points. Then the relative
Teichm\"uller space $\Teich(u) = {\cal J}(u)/\Equiv_0(u)$ is weakly
contractible (i.e.\ all its homotopy groups are trivial).
\end{statement}

\begin{proof}
Apparently, for any complex structure $I$ on $N$ there exists exactly one
$u$-compatible pair ${\cal I} \in \Compl(u)$ such that $I = p_2({\cal I})$.
This, together with the Corollary of Statement~\ref{St:E0Contr}, allows to
consider the space $\Teich(u)$ as a quotient of the space $\cal J$ of
complex structures at $N$ by action of the group $\Equiv_0\two(u)$.  The
point stabilizer of this action consists of diffeomorphisms of $N$
homotopic to identity, preserving branching values of $u$ and holomorphic
with respect to the relevant complex structure. Since $B(u)$ contains at
least three points, the stabilizer is trivial.

The natural quotient mapping ${\cal J} \to \Teich(u)$ is thus a fibration
with the fiber $\Equiv_0\two(u)$. The fiber is contractible by Statement
~\ref{St:E0Contr}, and the total space is contractible, too. An exact
homotopy sequence of the fibration ${\cal J} \to \Teich(u)$ shows now that
$\Teich(u)$ is weakly contractible.
\end{proof}

\begin{Remark}
In fact, $\Teich(u)$ is contractible. What we need, though, is only the
fact that it is simply connected (see the proof of Theorem
~\ref{Th:Pi1IsMap} below).
\end{Remark}

{\def \proofName {Proof of Theorem~\ref{Th:Pi1IsMap}}
\begin{proof}
follows immediately from Statement~\ref{St:TContr}: the Teichm\"uller
space $\Teich(u)$ is simply connected, and one has only to observe that the
action of the group $\Map(u)$ in $\Teich(u)$ is faithful and discrete.
\end{proof}
}

\section{Equivariant mapping class group}\label{Sec:Equiv}

From now on we restrict ourselves to the model example, a ramified covering
$u: M_g \to \Complex P^1$ described in Introduction. Recall that $q$ is a
homogeneous polynomial of two complex variables, of degree $2g+1$. The
associated polynomial of one variable $Q$ is defined by the equation
$q(x,z) = z^{2g+1} Q(x/z)$. We suppose that all the $2g$ critical values of
the polynomial $Q$ are distinct.

The full preimage $u^{-1}(a)$ of a generic point $a \in \Complex P^1$
consists of $4g+2$ points. There are $2g+2$ branching values of $u$ having
less preimages. Namely, $u^{-1}(\infty)$ contains only one point,
$u^{-1}(0)$ contains $2g+1$, and for each of the $2g$ critical values $P_0,
\dots, P_{2g-1}$ of the polynomial $Q$, the set $u^{-1}(P_i)$ has
cardinality $4g$. We suppose that $g > 0$, and thus mapping $u$ satisfies
hypotheses of Theorem~\ref{Th:Pi1IsMap}. In this Section we study the
structure of the group $\Map(u) = \pi_1(\Mod(u))$.

To avoid confusion we will always denote the target manifold of the mapping
$u$ as $\Complex P^1$. One should remember, though, that we do not assume a
complex structure on it to be fixed.

\subsection{Mapping classes of sphere with marked
points}\label{SSec:MapSphere}

Consider the projection $p_2: \Equiv(u) \to \Equiv\two(u)$. Mappings of
$\Complex P^1$ to itself that belong to the image of $p_2$ preserve the
branching order. So, they map points $0$ and $\infty$ to themselves, and
points $P_0, \dots, P_{2g-1}$ are permuted in some way. In other words,
$\Equiv\two(u)$ is a subgroup of $\Diff(\Complex P^1, \{0\}, \{\infty\},
\{P_0, \dots, P_{2g-1}\})$. The mapping class group $\Map(\Complex P^1,
\{0\}, \{\infty\}, \{P_0, \dots, P_{2g-1}\})$ will be called $G_{2g}$ for
short. We are now to study the group $G_{2g}$.

A simple curve $\phi$ joining the points $0$ and $\infty$ of $\Complex P^1$
and avoiding the points $P_0, \dots, P_{2g-1}$ will be called a {\em
meridian}. A system $\phi_0, \dots, \phi_{2g-1}$ of $2g$ meridians
(numbered) will be called a {\em slicing} if the meridians do not intersect
each other (except at $0$ and $\infty$), and cut $\Complex P^1$ into
``slices'' each containing exactly one point $P_k$. The group $G_{2g}$ acts
on the set of homotopy classes of slicings. It follows immediately from
Alexander's lemma that this action is faithful (a point stabilizer is
trivial).

Since we are interested with homotopy classes of mappings only, we may
suppose that $P_k = \exp((2k+1)\pi i /2g) \in \Complex P^1,\, k = 0, \dots,
2g-1$. We can also fix a ``standard'' slicing where the meridian $\phi_k$
is the line $\name{arg} z = 2\pi k/2g$. Thus, an element $\a \in G_{2g}$ is
completely determined by the homotopy classes of meridians $\a(\phi_0),
\dots, \a(\phi_{2g-1})$.

Consider a rotation $z \mapsto z \exp(2\pi i/2g)$. It shifts cyclically
meridians of the standard system: $\phi_0 \mapsto \phi_1 \mapsto \dots
\mapsto \phi_{2g-1} \mapsto \phi_0$. The corresponding mapping class will
be called $\lambda \in G_{2g}$. Obviously, $\lambda$ is an element of order
$2g$.\label{Mk:DefLambda}

Consider now a subgroup $G^0_{2g} \subset G_{2g}$ consisting of elements
mapping the meridian $\phi_0$ to itself.

\begin{statement}\label{St:GenBraid}
The group $G^0_{2g}$ is isomorphic to a braid group $B_{2g}$ on $2g$
strands.
\end{statement}

\begin{figure}
\centerline{\epsfbox{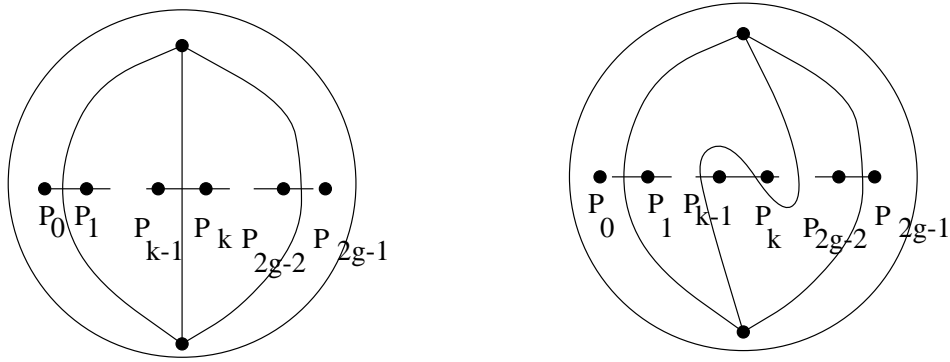}}
\caption{Mapping class $u_k$}\label{Fig:Gener}
\end{figure}

See \cite{Birman} for proof. A construction found in~\cite{PresBraid}
relates to every tree with the vertices $P_0, \dots, P_{2g-1}$ a system of
generators of the group $G^0_{2g}$. A ``line'' tree corresponds to the
standard generators $u_1, \dots, u_{2g-1}$ of the braid group; an action of
$u_k$ is shown in Fig.~\ref{Fig:Gener} (it is assumed that $u_k$ is
identical outside a small neighborhood of the segment $P_{k-1}P_k$).

\begin{statement}\label{St:Gener}
The mapping class group $G_{2g}$ is generated by its subgroups $G^0_{2g}$
and $\langle \lambda\rangle \subset G_{2g}$.
\end{statement}

See \cite{Cohen} for proof.

Let $s_k$ be a simple curve which starts at $\infty$, encircles the point
$P_k$, and then ends at $\infty$ not intersecting the meridians $\phi_i,\,
i = 0, \dots, 2g-1$. Homotopy classes of the loops $s_k$ are free
generators of the fundamental group $\pi_1(\Complex P^1 \setminus \{0, P_0,
\dots, P_{2g-1}\}; \infty)$ ($\infty$ is a base point). Group $G_{2g}$ acts
on this fundamental group in the following way:
\begin{eqnarray}
&&\lambda s_k = s_{k+1},\label{Eq:ActLambda}\\
&&u_k s_{k-1} = s_k,\label{Eq:ActU1}\\
&&u_k s_k = s_k^{-1} s_{k-1} s_k,\label{Eq:ActU2}\\
&&u_k s_l = s_l, \quad \mbox{if\ } l \ne k-1, k.\label{Eq:ActU3}
\end{eqnarray}

Addition in subscripts is made here $\name{mod}\,2g$.

Automorphisms $\lambda$ and $u_k,\, k = 1, \dots, 2g-1$ of the free group
$F(s_0, \dots, s_{2g-1})$ given by
formulas~(\ref{Eq:ActLambda})--(\ref{Eq:ActU3}) generate a subgroup of the
group $\name{Aut}(F(s_0, \dots, s_{2g-1}))$ called {\em braid-cyclic group}
$BC_{2g}$ on $2g$ strands. Thus
formulas~(\ref{Eq:ActLambda})--(\ref{Eq:ActU3}) define a homomorphism
$G_{2g} \to BC_{2g}$. It follows easily from Alexander's lemma that this
homomorphism has no kernel, i.e.\ is an isomorphism.

Thus we can summarize the results about group $G_{2g}$ in the following

\begin{theorem}\label{Th:BraidCycl}
The mapping class group $G_{2g} \bydef \Map(\Complex P^1; \{0\},
\{\infty\}, \{P_0, \dots, P_{2g-1}\})$ is isomorphic to the braid-cyclic
group $BC_{2g}$. It is a product (not direct) $G_{2g} = \langle
\lambda\rangle G^0_{2g}$. Here the subgroup $G^0_{2g}$ is a stabilizer of
the meridian $\phi_0$, it is isomorphic to the braid group $B_{2g}$. The
element $\lambda$ has order $2g$.
\end{theorem}

\subsection{Liftable mapping classes}\label{SSec:Lift}

In this Subsection we are to answer a question, which mapping classes $\a
\in G_{2g}$ are liftable, i.e.\ belong to the image of the homomorphism
$p_2: \Map(u) \to G_{2g}$. The subgroup of liftable classes will be
denoted $\Lift(u)$.

Consider the full preimage $u^{-1}(\phi) \subset M_g$ of a meridian $\phi$
under the mapping $u$. It is a graph having $2g+2$ vertices. One of them is
a point $T = [0:1:0] = u^{-1}(\infty)$ and the other $2g+1$ are elements
$S_1, \dots, S_{2g+1}$ of $u^{-1}(0)$ (zeros of the polynomial $Q$). The
graph has $4g+2$ edges, each one joining $T$ with some $S_k$. The valence
of the vertex $T$ is $4g+2$, and the valence of each vertex $S_k$ is $2$.
This allows us to ignore vertices $S_k$ and to consider $u^{-1}(\phi)$ as a
graph with a single vertex $T$ and $2g+1$ loops $\gamma_1, \dots,
\gamma_{2g+1}$ attached to it. Vertices $S_k$ look then as midpoints of the
loops. Relative positions of the loops are described in

\begin{statement}\label{St:Polygon}
If one cuts $M_g$ along the $2g$ curves $\gamma_1, \dots, \gamma_{2g}$, one
obtains a $4g$-gon embedded into $M_g$ so that its opposite sides are glued
together (correspond to the same curve $\gamma_i$). The remaining curve
$\gamma_{2g+1}$ joins opposite vertices of the $4g$-gon obtained (looks
like its ``large diagonal'').
\end{statement}

Indeed, numeration of the curves $\gamma_i$ can be chosen at random, and so
any $2g$ curves $\gamma_i$ can become $\gamma_1, \dots, \gamma_{2g}$
mentioned in the Statement.

\begin{proof}
Take the ``square root'' of the mapping $u$ presenting it as $u = s \circ
v$ where $v: M_g \to \Complex P^1$ maps the point $[x:y:z]$ to $[y:z]$, and
$s: \Complex P^1 \to \Complex P^1$ is defined as $s(\xi) \bydef \xi^2$.
Let $\Phi$ be the full preimage of the meridian $\phi$ under $s$. Without
loss of generality one can suppose that $\Phi$ is the real axis. It divides
$\Complex P^1$ into two open disks ${\cal U}_\pm = \{\xi \mid \name{Im} \xi
\mathrel{{<}{>}} 0\}$. Preimages of the disks ${\cal U}_\pm$ under the
mapping $v$ are the sets $D_\pm = v^{-1}({\cal U}_\pm) = \{(x,y) \in
\Complex^2 \mid y^2 = Q(x),\, \name{Im} y \mathrel{{<}{>}} 0\}$. Projection
$(x,y) \mapsto x$ maps homeomorphically each of the sets $D_\pm$ to the set
$A = \{x \in \Complex \mid Q(x) \notin \Real_+\}$. The polynomial $Q$ does
not have critical values on the real axis $\Phi$ ($\Phi$ is a preimage of
the meridian $\phi$), and so, $A$ is homeomorphic to a disk, and the same
is true for $D_\pm$, too. Thus we have proved that the $u^{-1}(\phi)$ cuts
$M_g$ into a union of two disks.

The curves $\gamma_1, \dots, \gamma_{2g+1}$ constitute the boundary
of each disk $D_\pm$. So, orientation of $D_+$ fixes a cyclic order of the
curves $\gamma_1, \dots, \gamma_{2g+1}$, and the orientation of $D_-$ does
the same thing. Consider the involution $\mu:  M_g \to M_g$ given by the
formula $\mu(x,y) = (x, -y)$. It exchanges the disks $D_+ \leftrightarrow
D_-$, leaves every curve $\gamma_i$ fixed (as a whole) and preserves the
orientation. Thus, the cyclic orderings induced in the set $\gamma_1,
\dots, \gamma_{2g+1}$ by orientations of the disks $D_+$ and $D_-$
are the same. This finishes the proof.
\end{proof}

The last assertion of Statement~\ref{St:Polygon} was proved using the
involution $\mu$. This involution looks like a central symmetry of the
$4g$-gon described in Statement~\ref{St:Polygon}. It is clear that $p_2
\circ \mu  = p_2$, so that $\mu$ projects to the identity mapping of
$\Complex P^1$. By this reason the full preimage $u^{-1}({\cal U}) \subset
M_g$ of any subset ${\cal U} \subset \Complex P^1$ is centrally
self-symmetric. So, to describe $u^{-1}({\cal U})$ it suffices to describe
its intersection with one disk, say, $D_+$.

Statement~\ref{St:Polygon} allows to view the disk $D_+$ as a polygon with
$2g+1$ vertices. Take now into account the $2g+1$ points $u^{-1}(0)$ and
regard them as vertices, too, to obtain a $(4g+2)$-gon. We supply the
vertices of this $(4g+2)$-gon with alternating coloring, so that the $2g+1$
vertices are $u^{-1}(0)$ black, and the $2g+1$ vertices $u^{-1}(\infty)$
are white. All the white vertices are indeed the same point in the $M_g$,
but Statement~\ref{St:Polygon} allows to regard them as different ones. We
will do it systematically, without further notice. For example, if two
smooth curves $\gamma_1, \gamma_2 \subset D_+$ intersect at the point
$u^{-1}(\infty)$ only, and are separated near it by some curve $\gamma_i
\subset u^{-1}(\phi)$, then we say that $\gamma_1$ and $\gamma_2$ have no
common points. They look as starting at different vertices of the polygon
$D_+$.

Statement~\ref{St:Polygon} describes the preimage of an individual
meridian. The following theorem deals with preimages of slicings.

\begin{theorem}\label{Th:PreimCompl}
Let $\phi_1, \dots, \phi_{2g}$ be a slicing. Let the preimage
$u^{-1}(\phi_1)$ cut $M_g$ into a union $D_+ \sqcup D_-$. Let $\Gamma
\bydef u^{-1}(\phi_2 \cup \dots \cup \phi_{2g}) \cap D_+$. Then $\Gamma$ is
a union of simple curves (edges) starting in white vertices and ending at
black ones. The edges have no common internal points. They cut $D_+$ into
parts (faces), each of them being homotopic either to a quadrangle or to a
lune (biangle). Each face contains exactly one preimage of some point $P_i
\in \Complex P^1$.
\end{theorem}

\begin{proof}
Recall that the meridians $\phi_0, \dots, \phi_{2g-1}$ of the slicing cut
$\Complex P^1$ into slices each one containing exactly one point $P_i$.
Connected components of preimages of these slices are the parts into
which $\Gamma$ divides $M_g$.

Consider, for example, the slice ${\cal U}$ bounded by $\phi_1$ and
$\phi_2$, and containing $P_1$. Split it into a union of curves
$\phi_\tau$, $1 \le \tau \le 2$, where the curves $\phi_\tau,\, \tau \ne
3/2$ are meridians having no common internal points, and the curve
$\phi_{3/2}$ starts at $\infty$, passes $P_1$, ends at $0$, and has no
common internal points with the other curves $\phi_\tau$ either. Preimages

$u^{-1}(\phi_\tau) \cap D_+$ for $\tau \ne 3/2$ are systems of smooth
curves joining black and white vertices of $D_+$. These curves are
homotopic to the respective components of $u^{-1}(\phi_1)$ (for $1 \le \tau
< 3/2$) or $u^{-1}(\phi_2)$ (for $3/2 < \tau \le 2$).

Preimage of the point $P_1$ in the disk $D_+$ consists of $2g$ points. One
of them is a branching point of the multiplicity $2$, and the other $2g-1$
points are not critical. Thus the preimage $u^{-1}(\phi_{3/2})$ consists of
$2g$ components. $2g-1$ of them are smooth, and the corresponding part of
$M_g$ is a lune. The last component looks like two smooth curves
intersecting transversally at one point. Apparently, the corresponding part
of $M_g$ is a quadrangle.
\end{proof}

Recall that the meridians of the slicing are numbered: $\phi_0, \phi_1,
\dots, \phi_{2g-1}$. Number the points $P_i$ also, so that the slice bounded
by $\phi_k$ and $\phi_{(k+1) \bmod 2g}$ contain $P_k$. The last assertion
of Theorem~\ref{Th:PreimCompl} allows then to label the faces with symbols
$0, 1, \dots, 2g-1$ --- the face $\cal F$ containing the preimage of $P_i$
is labelled $i \bydef \ell({\cal F})$.

Let us direct each edge, as well as each side of the $(4g+2)$-gon, from its
white end to its black end.  Let ${\cal F}_1$ and ${\cal F}_2$ be two faces
separated by an edge $e$, ${\cal F}_1$ lying to the right of it, and ${\cal
F}_2$, to the left. Since the mapping $u$ preserves the orientation, labels
of the faces satisfy the equation $\ell({\cal F}_2) =  \ell({\cal F}_1) + 1
\pmod{2g}$. This remark allows to reproduce positions and labels of the
lune-like (biangular) faces uniquely up to a homotopy, if positions and
labels of the quadrangular faces are given. Thus, the full preimage of a
slicing can be described, up to a homotopy, by two sorts of data: an
embedding of the $(2g+1)$-loop graph $u^{-1}(\phi_1)$ into $M_g$, and a
splitting of the corresponding $(4g+2)$-gon $D_+$ into quadrangles labelled
$0, 1, \dots, 2g-1$. Recall, too, that the vertices of the $(4g+2)$-gon
$D_+$ are colored alternatingly black and white. Splitting of $D_+$ into
labelled quadrangles together with vertex coloring will be called a
(marked) {\em quadrangulation}.  Quadrangulation corresponding to the
slicing $\Phi$ will be denoted $\square(\Phi)$.

Fix now a ``standard'' slicing $\Phi_0 = \{\phi_0^0, \dots,
\phi_{2g-1}^0\}$, for example, the one described in Subsection
~\ref{SSec:MapSphere}. Recall that the mapping class group $G_{2g} =
\Map(\Complex P^1, \{0\}, \{\infty\}, \{P_0, \dots, P_{2g-1}\})$ acts
(faithfully) at the set of homotopy classes of slicings. Now we can
formulate the required criterion of liftability:

\begin{theorem}\label{Th:LiftStab}
The element $\a \in G_{2g}$ is liftable (belongs to the subgroup $\Lift(u)
\subset G$) if and only if $\square(\a(\Phi_0)) = \square(\Phi_0)$.
\end{theorem}

\begin{Proof}
follows easily from Alexander's lemma.
\end{Proof}

\subsection{A short exact sequence}\label{SSec:ExaSeq}

\begin{Lemma}\label{Mk:Autom}
Let $Q:\Complex \to \Complex$ be a polynomial of odd degree $2g+1$ with
$2g$ distinct critical values, and $f: \Complex \to \Complex$ be a
diffeomorphism such that $Q(f(x)) = Q(x)$ for all $x \in \Complex$. Then
$f(x) \equiv x$.
\end{Lemma}

\begin{proof}
The mapping $f$ in question is an algebraic diffeomorphism of $\Complex$,
and thus a linear function: $f(x) = px + q$. This function preserves all
the $2g$ critical points (distinct) of the polynomial $P$, and therefore
$f(x) \equiv x$.
\end{proof}

Recall (see Theorem~\ref{Th:BraidCycl}) that the group  $G_{2g}$ is an
indirect product $\langle \lambda \rangle G^0_{2g}$. Here the cyclic
subgroup $\langle \lambda \rangle$ generated by the ``rotation'' $\lambda$
(see definition at page~\pageref{Mk:DefLambda}) of order $2g$. The subgroup
$G^0_{2g}$ is isomorphic to the braid group $B_{2g}$ and consists of the
mapping classes preserving the meridian $\phi_0$ of the standard slicing.
Denote $\Lift^0(u) = \Lift(u) \cap G^0_{2g}$.

\begin{theorem}\label{Th:ExaSeq}
The group $\Map(u)$ fits into a short exact sequence of groups
\begin{equation}\label{Eq:Exact}
0 \to \Integer_2 \to \Map(u) \to \Lift(u) \to 0
\end{equation}
The generator of the group $\Integer_2$ is mapped to the involution $\mu$
lying in the center of $\Map(u)$. The second arrow is the projection $p_2$.
The group $\Map(u)$ is an indirect product $R F$ where $p_2(F) =
\Lift^0(u)$, and $p_2: F \to \Lift^0(u)$ is an isomorphism. The group $R$
is a cyclic group generated by the element $(\Lambda, \lambda) \in \Map(u)$
of order $4g$.
\end{theorem}

Recall that the involution $\mu: M_g \to M_g$ acts by the formula $\mu(x,y)
= (x,-y)$ (see remark after the proof of Statement~\ref{St:Polygon}).

\begin{proof}
Prove first that the kernel of the mapping $p_2: \Map(u) \to \Lift(u)$ is
generated by the involution $\mu$. Projections $p_1(\tau)$ and $p_2(\tau)$
determine an element $\tau \in \Map(u)$ completely. Thus, it suffices to
prove that $\mu$ is the only nontrivial element of the group
$\Equiv\one(u)$ which projects to the identity mapping of $\Complex P^1$.
Apparently any such mapping $\tau: M_g \to M_g$ acts by the following way:
$\tau(x,y) = (f(x), \pm y)$ where $f: \Complex \to \Complex$ is a
diffeomorphism preserving the polynomial $Q$: $Q(f(x)) = Q(x)$. By the
lemma, $f(x) \equiv x$.

Fix now a standard slicing $\Phi_0$ at $\Complex P^1$, and let $M_g = D_+
\cup D_-$ like in Theorem~\ref{Th:PreimCompl}. Define $F$ as the following
subgroup of $\Map(u)$: $f \in F$ $\Leftrightarrow$ $p_1(f)$ preserves $D_+$
(and therefore $D_-$, too). Apparently, $p_2(F) = \Lift^0(u)$. Since $\mu
\notin F$, the restriction $\left.p_2\right|_F$ is an isomorphism.

Let now $(\Lambda, \lambda)$ be a lifting of the element $\lambda \in
G_{2g}$. Such lifting exists (i.e.\ $\lambda \in \Lift(u)$): the mapping
$\Lambda:M_g \to M_g$ is a rotation of the $4g$-gon $D_+ \cup D_-$ to the
$1/4g$ of the full cycle (more exactly, it is one of two mappings
projecting to $\lambda$, the second in $\Lambda \mu$). Obviously, $\Lambda$
is an element of order $4g$. Theorem~\ref{Th:BraidCycl} implies that the
group $\Lift(u)$ is generated by the element $\lambda = p_2((\Lambda,
\lambda))$ and the subgroup $\Lift^0(u) = p_2(F)$; the kernel of $p_2$ is
generated by $(\mu, \name{Id}) = (\Lambda, \lambda)^{2g}$. So, the group
$\Map(u)$ is generated by $(\Lambda, \lambda)$ and the subgroup $F$.
\end{proof}

\section{Actions of the braid-cyclic group}\label{Sec:Action}

\subsection{Action on marked quadrangulations}\label{SSec:ActQuad}

The liftability criterion, as formulated in Theorem~\ref{Th:LiftStab}, does
not answer directly a question when an element $\a \in G_{2g} = BC_{2g}$
written in generators $\lambda, u_1, \dots, u_{2g-1}$ is liftable. In this
Section we will try to do it giving an algorithmic description of the
action of $\alpha$ at marked quadrangulations.

Consider a quadrangulation $\Delta$ of the $(2N+2)$-gon $S$ with vertices
colored alternatingly black and white. It is easy to see that any
quadrangle has two black and two white vertices, and that its opposite
vertices have the same color. Fix now a vertex $v$. Orientation of $S$
determines a linear ordering of the quadrangles (faces of $\Delta$) having
$v$ their vertex. On the other hand, this set is ordered by the labels
$0,1, \dots, N-1$ of faces. We call quadrangulation $\Delta$ {\em monotone}
if these two orderings coincide for all white vertices $v$ and are opposite
for all black vertices.

\begin{statement}\label{St:QuadMono}
For any slicing $\Phi$ the quadrangulation $\square(\Phi)$ of the
$(4g+2)$-gon $S = D_+$ is monotone.
\end{statement}

\begin{proof}
Consider the quadrangulation $\square(\Phi)$ together with the biangular
faces. Mapping $u: M_g \to \Complex P^1$ preserves the orientation, and
therefore the cyclic ordering of edges in every white vertex of the graph
$\Gamma = u^{-1}(\Phi)$ coincides with the ordering of the meridians:
$\phi_0, \phi_1, \dots, \phi_{2g-1}$. So, if the quadrangulation is not
monotone then there exists a pair of faces (biangular or quadrangular)
attached to the same vertex and separated by the preimage of the meridian
$\phi_0$. This is impossible because $u^{-1}(\phi_0)$ forms the boundary
of the $(4g+2)$-gon $S$. The same argument applies to black vertices.
\end{proof}

Note now that for any monotone quadrangulation the face labelled $1$ has
only one or two sides adjacent to other faces, and if there are two such
sides then they are opposite. The other two or three sides lie in the
boundary of the polygon $S$. The same remark applies to the face labelled
$N$. Now we are ready to define an action of the braid-cyclic group $BC_N$
on the set of monotone quadrangulations.

\begin{statement}\label{St:ActQuadr}
The following rules define an action of the group $BC_N$ on the set of
quadrangulations of the $(2N+2)$-gon, with faces labelled $0, 1, \dots,
N-1$.
\begin{enumerate}
\item Element $u_k$ acts only on the faces labelled $k-1$ and $k$. Here are
two cases:
\begin{enumerate}
\item If the faces labelled $k-1$ and $k$ are not adjacent then they
exchange their labels (the quadrangulation remains the same).
\item If these faces are adjacent, and thus form a hexagon, then this
hexagon is rotated by $1/3$ of the full cycle.
\end{enumerate}
\item Element $\lambda$ acts as follows. First, it shifts the face
labelling cyclically: $1 \mapsto 2 \mapsto \dots \mapsto N \mapsto 1$.
Let $\cal F$ be the face labelled $1$, and $AB$ and $CD$ be its
opposite side lying on the boundary of $S$ (they exist by the previous
remark), and let $P$ and $Q$ are vertices of $S$ such that $A, B, P$ are
subsequent vertices, and $C, D, Q$, too. Then the face ${\cal F} = ABCD$ is
replaced by the face $PBQD$.
\end{enumerate}
\end{statement}

\begin{figure}
\centerline{\epsfbox{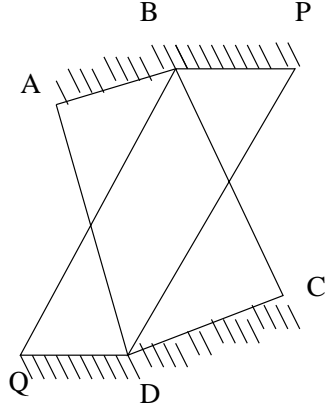}}
\caption{Action of $\lambda$ at a quadrangulation}\label{Fig:ActLAmbda}
\end{figure}

Action of $\lambda$ on quadrangulations (``the flip'') is shown in
Fig.~\ref{Fig:ActLAmbda}. Action of $u_k$ is shown at Fig.~\ref{Fig:Rotat}
below.

\begin{proof}
To prove that we have defined an action, we should check the defining
relations of the group $BC_N$ (see~\cite{Birman} for proof; it follows
easily from relations~(\ref{Eq:ActLambda})--(\ref{Eq:ActU3}), too):
\begin{eqnarray}
&&\lambda^N = 1\label{Eq:LambdaN}\\
&&\lambda u_k = u_{k+1}\lambda, \quad \quad k = 1, \dots,
N-2.\label{Eq:RelLamU}\\
&&u_k u_l = u_l u_k, \quad k,l = 1, \dots, N-1, \lmod k-l\rmod \ge
2;\label{Eq:RelBrComm}\\
&&u_k u_{k+1} u_k = u_{k+1} u_k u_{k+1}, \quad k = 1, \dots,
N-2.\label{Eq:RelABA}
\end{eqnarray}
The checking is routine.
\end{proof}

Let $a \mapsto a^{inv}$ be an anti-homomorphous involution of the group
$BC_N$ sending each $u_k$ to itself, and $\lambda$ to $\lambda^{-1}$.

\begin{theorem}\label{Th:TwoActions}
Actions of the group $BC_{2g}$ on the set of slicings and on the set of
monotone quadrangulations are connected by the following formula:
\begin{displaymath}
\square(\a(\Phi_0)) = \a^{inv}(\square(\Phi_0)).
\end{displaymath}
\end{theorem}

\begin{proof}
Write $\a = s_1 s_2 \dots s_M$ where each $s_i$ is either $\lambda$ or
$u_k^{\pm 1}$, and denote $\a' = s_2 \dots s_M$. To prove the theorem by
induction on $M$ suppose that for $\a'$ the theorem is proved, and consider
three cases. The symbol ${\cal F}_i$ denotes the face labelled $i$ in
the quadrangulation $\square(\Phi)$ where $\Phi = \alpha'(\Phi_0)$.

\begin{stages}{Case}
\nextstage $s_1 = u_k^{\pm 1}$, and the faces ${\cal F}_k$ and ${\cal
F}_{k+1}$ are not adjacent. When $u_k^{\pm 1}$ acts on $\Phi_0$, then all
the meridians except $\phi_k$ remain the same, points $P_k$ and $P_{k+1}$
exchange places, and each point $P_i$ goes to itself. Therefore the
quadrangulations $\square(\a'(\Phi_0))$ and $\square(\a(\Phi_0))$ can
differ in the faces labelled $k$ and $k+1$ only. These faces are not
adjacent, so the quadrangulation remains geometrically the same. Labels $k$
and $k+1$ exchange places because the points $P_k$ and $P_{k+1}$ do so.

\nextstage $s_1 = u_k$, and the faces ${\cal F}_k$ and ${\cal F}_{k+1}$ are
adjacent, thus forming a hexagon $H$. Like in the previous case, one shows
that the quadrangulations $\square(\a'(\Phi_0))$ and $\square(\a(\Phi_0))$
can differ inside the hexagon $H$ only. Denote $\rho'$ the border
between ${\cal F}_k$ and ${\cal F}_{k+1}$. The image $u(\rho')$ is the
meridian $\a'(\phi_k)$. Image of some other diagonal of $H$, say, $\rho$,
is the meridian $\a(\phi_k)$.

\begin{figure}
\centerline{\epsfbox{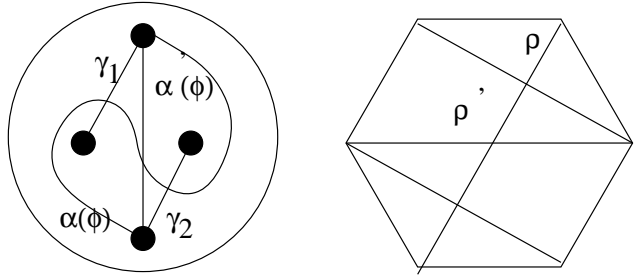}}
\caption{Rotation of a hexagon}\label{Fig:Rotat}
\end{figure}

Draw curves $\gamma_1$ and $\gamma_2$ like in Fig.~\ref{Fig:Rotat}. The
sets $\ell_1 = u^{-1}(\gamma_1) \cap H$ and $\ell_2 = u^{-1}(\gamma_2) \cap
H$ are homotopic to diagonals of ${\cal F}_1$ and ${\cal F}_2$,
respectively. Then, diagonal $\rho$ has a common point with the curve
$\ell_1$, and another one, with $\ell_2$, and these points are preserved
under any homotopy of curves $\ell_1$ and $\ell_2$ with fixed endpoints. It
means that $\rho$ is the diagonal obtained from $\rho'$ by a
counterclockwise rotation of $H$ to the $1/3$ of the full cycle. To
determine the labelling of the quadrangles obtained we can use the same
arguments as in the previous case.

\nextstage $s_1 = \lambda$. In this case the meridians are mapped
cyclically: $\phi_0 \mapsto \phi_1 \mapsto \dots \mapsto \phi_{2g} \mapsto
\phi_0$. Thus the splitting of the whole $M_g$ into quadrangles does not
change, but the face which was labelled $i$ is now labelled $i-1$. Besides
this, the splitting $M_g = D_+ \cup D_-$ changes: we cut $M_g$ along the
preimage of the new $\phi_0$ which was previously $\phi_1$. One can easily
observe that this corresponds exactly to the ``flip'' of quadrangulation
defining the action of $\lambda^{-1}$ in Statement~\ref{St:ActQuadr}.
\end{stages}
\end{proof}

\begin{Remark}
Indeed, we can choose any slicing $\Phi$ as a ``standard'' slicing
$\Phi_0$. Another choice of $\Phi_0$, though, leads to another system of
generators $\lambda, u_1, \dots, u_{2g-1}$, and thus to another involution
$a \mapsto a^{inv}$, so that Theorem~\ref{Th:TwoActions} remains true. If
one changes $\Phi_0$ only (but involution remains the same), the theorem
fails. See \cite{PresBraid} for more information on systems of generators
of the braid group.
\end{Remark}

\subsection{Action on trees with labelled edges}\label{SSec:ActTree}

In this Subsection we define an action of the group $BC_N$ on the set of
all trees with $N$ edges labelled $0, 1, \dots, N-1$, and investigate its
relation with the action of Statement~\ref{St:ActQuadr}. This action is of
some independent interest, and also allows us to obtain some information
about the group $\Lift(u)$.

\begin{statement}\label{St:ActTree}
The following rules define an action of the group $BC_N$ on the set of all
trees with $N$ edges labelled $0, 1, \dots, N-1$.

\begin{enumerate}
\item Element $u_k$ acts only on the edges labelled $k-1$ and $k$. Here are
two cases:
\begin{enumerate}
\item If the edges labelled $k-1$ and $k$ do not have common vertices then
they exchange their labels (the tree remains the same).
\item If the edge labelled $k-1$ joins the vertices $A$ and $B$, and the
edge labelled $k$ joins $A$ and $C$ then the edge $AB$ receives label $k$,
edge $AC$ is erased, and the edge $BC$ is drawn and labelled $k-1$.
\end{enumerate}
\item Element $\lambda$ does not affect the tree but shifts the edge
labelling cyclically: $1 \mapsto 2 \mapsto \dots \mapsto N \mapsto 1$.
\end{enumerate}
\end{statement}

The proof copies that of Statement~\ref{St:ActQuadr}.

\begin{statement}\label{St:Trans}
Action of the braid-cyclic group $BC_N$ on the set of trees with edges
labelled $0, 1, \dots, N-1$ is transitive.
\end{statement}

\begin{proof}
We prove that already the action of the subgroup $B_N \subset BC_N$
generated by $u_1, \dots, u_{N-1}$ is transitive.

Take some vertex $a$ of the tree as a root, and orient every edge away from
the root (downwards). We can speak now about upper and lower end of an
edge. Two edges having the same upper end will be called {\em brothers}.
If the lower end of one edge is the upper end of another edge, then such
edges will be called {\em parent} and {\em child}, respectively.

Call the {\em complexity} of a tree $\Gamma$ (notation $C(\Gamma)$) the sum
of lengths of all strictly decreasing paths starting at the root. The
minimal possible value for the complexity is $N$, the number of edges. Only
one tree has complexity $N$, this is the ``bush'' tree where all $N$ edges
are attached to the root. We will prove that for any other tree $\Gamma$
there exists an element $\a \in B_N$ such that $C(\a(\Gamma)) < C(\Gamma)$.
This will mean that $\Gamma$ can be mapped to the ``bush'' tree by some
element of the braid group, and thus the action is transitive.

The generator $u_k$ of the braid group affects the tree complexity in the
following way (the edge labelled $i$ is denoted $e_i$):
\begin{enumerate}
\item\label{It:NonAdj} If $e_{k-1}$ and $e_k$ are not adjacent then
$C(u_k(\Gamma)) = C(\Gamma)$.
\item\label{It:Broth} If $e_{k-1}$ and $e_k$ are brothers then
$C(u_k(\Gamma)) > C(\Gamma)$.
\item\label{It:ParMinus} If $e_{k-1}$ is the parent of $e_k$ then
$C(u_k(\Gamma)) < C(\Gamma)$.
\item\label{It:ChiMinus} If $e_k$ is the parent of $e_{k-1}$ then any
inequality between $C(u_k(\Gamma))$ and $C(\Gamma)$ is possible, but
$C(u_k^2(\Gamma)) = C(u_k^{-1}(\Gamma)) < C(\Gamma)$.
\end{enumerate}
We will look for such element $\a = u_{i_1} \dots u_{i_M}$ that for $s = 1,
\dots, M-1$ action of $u_{i_s}$ does not change complexity, and action of
$u_{i_M}$ reduces it. To find $\a$ consider several cases.

\begin{enumerate}
\item\label{It:NoN} There exists an edge $e_k$ with $k < N$ and having a
child.
{\def \theenumii {\arabic{enumi}.\arabic{enumii}}
\def \labelenumii {\theenumii}
\expandafter\def\csname p@enumii\endcsname {}
\begin{enumerate}
\item\label{It:DecrPoss} $e_{k+1}$ is a child of $e_k$. In this case we
apply $u_{k+1}$ and decrease complexity.
\item\label{It:ExiLess} $e_{k+1}$ is not a child of $e_k$ but there is a
child $e_s$ with $s < k+1$. Take the greatest such $s$. The edge $e_{s+1}$
is neither a brother nor a parent of $e_s$, so that we can apply $u_s$ not
increasing the complexity. Thus we have reproduced the same situation but
with $s \mapsto s+1$. Repeating the process several times we arrive to the
situation of Case~\ref{It:DecrPoss}.
\item\label{It:AllMore} All the children of $e_k$ have labels greater than
$k+1$. Then take the smallest label $s$. Here are two more cases:
{\def \theenumiii {\theenumii.\arabic{enumiii}}
\def \labelenumiii {\theenumiii}
\expandafter\def\csname p@enumiii\endcsname {}
\begin{enumerate}
\item\label{It:NotChi} $e_{s-1}$ is not a child of $e_s$. Then we apply
$u_s$ not changing the complexity. Thus we have reproduced the situation of
Case~\ref{It:AllMore} but with $s \mapsto s-1$. Repeating the process
several times we arrive to the situation of Case~\ref{It:DecrPoss}.
\item\label{It:Chi} $e_{s-1}$ is not a child of $e_s$. Here we apply
$u_s^2$ decreasing the complexity.
\end{enumerate}
}
\end{enumerate}
}
\item\label{It:N} Only the edge $e_N$ has children (the ``shrub'' tree).
Here are two cases:
{\def \theenumii {\arabic{enumi}.\arabic{enumii}}
\def \labelenumii {\theenumii}
\expandafter\def\csname p@enumii\endcsname {}
\begin{enumerate}
\item\label{It:N-1Chi} $e_{N-1}$ is a child of $e_N$. Here we can apply
$u_N^2$ to decrease complexity.
\item\label{It:N-1NotChi} $e_{N-1}$ is attached to the root and is not thus
a child of $e_N$. Let $e_k$ be a child of $e_N$ with the greatest $k$. So,
$e_{k+1}$ is neither a child nor a parent nor a brother of $e_k$, and we
can apply $u_{k+1}$ not changing the complexity. Thus we have reproduced
the same situation but with $k \mapsto k+1$. Repeating the process
several times we arrive to the situation of Case~\ref{It:N-1Chi}.\qed
\end{enumerate}
}
\end{enumerate}
\end{proof}

Establish now a link between the action of the braid-cyclic group on trees
and the action of Statement~\ref{St:ActQuadr}.

Let now $\Delta$ be a monotone marked quadrangulation of the $(2N+2)$-gon
with vertices colored alternatingly black and white. Recall that any
quadrangle has two black and two white vertices, its opposite vertices
having the same color. Draw a graph $\Gamma(\Delta)$ joining two black
vertices of every quadrangle. Edges of $\Gamma(\Delta)$ are naturally
labelled $0, 1, \dots, N-1$.

\begin{statement}\label{St:TreeToQuad}
The graph $\Gamma(\Delta)$ is a tree. For any tree $T$ with edges labelled
$0, 1, \dots, N-1$ there exists exactly one monotone quadrangulation
$\Delta$ such that $\Gamma(\Delta) = T$. Actions of the braid-cyclic group
on monotone quadrangulations and on trees are compatible:
$\Gamma(\a(\Delta)) = \a(\Gamma(\Delta))$ for any $\a \in BC_N$ and any
monotone quadrangulation $\Delta$.
\end{statement}

\begin{proof}
Suppose that $\Gamma$ has a cycle formed by edges $d_1, \dots, d_k$. Since
$d_1, \dots, d_k$ do not intersect, they form a polygon $P$. At least one
vertex of the original $(2N+2)$-gon lies inside $P$ which is impossible.
Thus, $\Gamma$ is a tree.

The second assertion is proved by induction on $N$. Consider the face $\cal
F$ of $\Delta$ marked $N-1$. The quadrangulation $\Delta$ is monotone, and
therefore two opposite sides of $\cal F$ are sides of the $(2N+2)$-gon.
Take now the tree $T$. The edge $e_{N-1}$ divides it into two subtrees, the
first containing $k$ edges, and the second, $N-k$ edges. This allows to
determine uniquely, up to a rotation of the $(2N+2)$-gon, the position
$\cal F$ --- its black diagonal divides the polygon into parts containing
$2k+3$ and $2N-2k+1$ vertices, respectively.  Positions of other faces are
then fixed uniquely, too, by the induction hypothesis.

The third assertion is proved by the obvious comparison of definitions, see
Statement~\ref{St:ActQuadr} and Statement~\ref{St:ActTree}.
\end{proof}

\begin{Corollary}
$\Lift(u)$ is a subgroup of the braid-cyclic group of finite index
$(2g+1)^{2g-2}$.
\end{Corollary}

\begin{proof}
Indeed, there exist $(n+1)^{n-2}$ different trees with $n$ numbered edges.
\end{proof}

For $N = 2g$ the graph $\Gamma(\Delta)$ can be also obtained by the
following construction. Consider a mapping $v: M_g \to \Complex P^1$
defined by the formula $v([x:y:z]) = [x:z]$. This mapping is not defined in
the point $[0:1:0]$ but outside it there is an equality $u = Q \circ v$
where $Q: \Complex P^1 \to \Complex P^1$ sends $\xi$ to $Q(\xi)$, $Q$ being
a polynomial used in definition of $M_g$. Consider a slicing $\Phi$ such
that $\Delta = u^{-1}(\Phi)$, and join every point $P_k$ with $0$ by a
simple curve $\gamma_k$ inside the corresponding slice. The preimage
$Q^{-1}\left(\bigcup_{k=0}^{2g-1} \gamma_k\right) \subset \Complex P^1$ was
considered in~\cite{Looij,Lyash,Arn}. It looks like a graph with two types
of edges. Edges of the first type join different preimages of the point
$0$.  The midpoints of these edges are preimages of the points $P_i$, one
edge for each $i$. Edges of the second type are ``hanging'' edges joining
a preimage of $0$ with a preimage of some $P_i$ (different from the
midpoints considered earlier). It was shown in~\cite{Arn} that the edges of
the first type form a tree $\cal G$. It is easy to see that $\Gamma(\Delta)
= v^{-1}({\cal G})$. This gives another proof of the first assertion of
Statement~\ref{St:TreeToQuad}.

Statement~\ref{St:Trans} implies that there exists a slicing $\Phi$ such
that the $\square(\Phi)$ is a ``trivial'' quadrangulation, where all the
quadrangles have a common black vertex. Consider such slicing and the
corresponding generators $\lambda, u_1, \dots, u_{2g-1}$ of the group
$G_{2g}$. It follows from Theorem~\ref{Th:LiftStab},
Statement~\ref{St:ActQuadr} and Theorem~\ref{Th:TwoActions} that the
element $U = u_1 \dots u_{2g-1}$ is liftable.

{\def \statementName {Conjecture}
\begin{Statement}
Elements $\lambda$ and $U$ generate the group $\Lift(u)$. There are no
relations between them except $\lambda^{2g} = 1$, and thus $\Lift(u)$ is
isomorphic to the free product $\Integer_{2g} * \Integer$.
\end{Statement}
}

\subsection{Action on subgroups of the free group}\label{SSec:Subgr}

In this Subsection we are to reveal the algebraic origin of the actions of
Statement~\ref{St:ActQuadr} and Statement~\ref{St:ActTree}. We show that it
is a part of the general action of the group $\name{Aut}(F_N)$ of
automorphisms of the free group $F_N$ on the set of its subgroups.
Embedding $BC_N \hookrightarrow \name{Aut}(F_N)$ is given by the
formulas~(\ref{Eq:ActU1})--(\ref{Eq:ActU3}). To establish the link we
identify first certain subgroups of $F_N$ with trees with labelled edges.

Consider $F_N$ as the fundamental group of the union of $N$ circles ${\cal
U} = \bigvee_{k=0}^{N-1} (S^1)_k$. Coverings spaces (nonramified)
$\Upsilon$ of ${\cal U}$ are graphs where all vertices have valence $2N$,
and edges are labelled $0, 1, \dots, N-1$. The label $k$ means that the
corresponding edge is mapped to the $k$-th circle of the union $\cal U$ by
the covering map $r:\Upsilon \to {\cal U}$. The number of vertices of
$\Upsilon$ is equal to the number of sheets of the covering $r$. If the
circles of the union $\cal U$ are oriented, then the graph $\Upsilon$
receives an orientation, too.

Recall that a {\em path} in the graph $\Upsilon$ is a sequence of edges
(two successive edges in the sequence have a common vertex) passing every
edge not more than once; {\em circuit} is a closed sequence of edges
passing every vertex not more than once.  Covering $\Upsilon$ will be
called {\em tree-like} if it possesses the following properties:
\begin{enumerate}
\item\label{It:1or2} All circuits in it have length $1$ or $2$.
\item\label{It:LabSame} If a circuit has length $2$ then both edges in it
have the same label.
\item\label{It:AllLab} For any number $k \in \{0, 1, \dots, N-1\}$ there
exists exactly one circuit of length $2$ whose edges are labelled $k$.
\end{enumerate}

For a tree-like covering space $\Upsilon$ a graph $\Gamma(\Upsilon)$ is
defined as follows. Vertices of $\Gamma(\Upsilon)$ are vertices of
$\Upsilon$. Two vertices, $A$ and $B$ of $\Gamma(\Upsilon)$ are connected
by the edge labelled $k$ if $\Upsilon$ contains edges $AB$ and $BA$ both
labelled $k$ (labels must coincide by Property~\ref{It:LabSame}). Loops
of $\Upsilon$ are ignored. By Property~\ref{It:1or2} $\Gamma(\Upsilon)$ is
a tree, and by Property~\ref{It:AllLab} it has $N$ edges so that every
label $0, 1, \dots, N-1$ is used exactly once.

The graph $\Upsilon$ can be restored uniquely by $\Gamma(\Upsilon)$.
Indeed, positions and labels of loops are determined unambiguously by
$\Gamma(\Upsilon)$ and the condition that for every vertex $v$ of
$\Upsilon$ and every $k \in \{0, 1, \dots, N-1\}$ exactly two edge ends
incident to $v$ should be labelled $k$.

Choose a base point $B$ in the covering space $\Upsilon$; without loss of
generality let $B$ be a vertex. Now we can relate the covering $r:\Upsilon
\to {\cal U}$ to the subgroup $\Group(\Upsilon) \bydef
r_*(\pi_1(\Upsilon,B)) \subset \pi_1({\cal U}, O) = F_N$ of the free group;
here $O$ is the vertex of the union. This correspondence between subgroups
and coverings with base point is one-to-one, by the following classical
theorem of homotopic topology (see e.g.~\cite{FuFo} for proof):

\begin{theorem}\label{Th:SubgrCov}
For any subgroup $G \subset F_N$ there exists covering $r:\Upsilon \to
{\cal U}$ with the base point $B$ such that $\Group(\Upsilon) = G$.  Any
two coverings $r_1, r_2$ possessing this property are equivalent in the
following sense: there exist a homeomorphism $f: \Upsilon_1 \to \Upsilon_2$
such that $f(B_1) = B_2$ and $r_2 \circ f = r_1$. The number of sheets of
the covering $r$ is equal to the index $[F_N:G]$ of the subgroup $G$.
\end{theorem}

Generators of the fundamental group $\pi_1({\cal U},O) = F_N$ are in
one-to-one correspondence with the circles of the union. This allows us to
think that the edges of the graphs $\Upsilon$ and $\Gamma(\Upsilon)$ are
labelled by generators $e_0, \dots, e_{N-1}$ of the group $F_N$, so that
it is possible to speak about products of labels, etc. By
Property~\ref{It:AllLab} for every $k$ there is exactly one edge of
$\Gamma(\Upsilon)$ labelled $e_k$, and we will call it simply ``edge
$e_k$''.

For any sequence of edges $\nu$ in the graph $\Upsilon$ define the element
$\tau_\nu \in F_N$ as $\tau_\nu = r_*(\nu)$. In other words, one should
multiply the labels of the edges passed by $\nu$. The label should be taken
with the exponent $+1$ if $\nu$ passes the edge in positive direction, and
with the exponent $-1$ (the inverse element) otherwise. The group
$\Group(\Upsilon)$ consists of all the elements $\tau_\nu$ for all closed
sequences of edges $\nu$ starting and ending at the base vertex $B$. We
will also use the notation $\tau_\nu$ where $\nu$ is a sequence of edges of
the graph $\Gamma(\Upsilon)$; here all the exponents are taken $+1$ since
the edges of $\Gamma(\Upsilon)$ bear no orientation.

It is easy to see that the following elements constitute a system of
generators for the group $\Group(\Upsilon)$:

\begin{enumerate}
\item Elements $a_{\nu,e}^{\Gamma(\Upsilon)} \bydef \tau_\nu e
\tau_\nu^{-1}$ where $\nu$ is an arbitrary path in $\Gamma(\Upsilon)$
starting at the base vertex $B$, and $e$ is a generator different from the
last letter in the word $\tau_\nu$.
\item Elements $b_{\nu,e}^{\Gamma(\Upsilon)} \bydef \tau_\nu e^2
\tau_\nu^{-1}$ where $\nu$ is an arbitrary path in $\Gamma(\Upsilon)$
starting at the base vertex $B$, not passing the edge labelled $e$, and
ending in the vertex $A$ attached to that edge.
\end{enumerate}

These elements are just $\tau_\phi$ for the following sequences of edges
$\phi$ in the graph $\Upsilon$:
\begin{enumerate}
\item For $a_{\nu,e}^{\Gamma(\Upsilon)}$ the sequence $\phi$ goes first
along the path $\nu$ in $\Gamma(\Upsilon)$. Every two successive vertices
of the path $\nu$ are connected in $\Upsilon$ by two edges having the same
label (by Property~\ref{It:LabSame}). Apparently, these edges must have
opposite orientations. The path $\phi$ takes each time the edge with the
positive direction. After this, $\phi$ passes the loop labelled $e$ and
attached to the final vertex of the path $\nu$, and then returns along the
path $\nu$ taking the same edges as for the first time.

\item For $b_{\nu,e}^{\Gamma(\Upsilon)}$ the sequence $\phi$ goes along the
path $\nu$ in $\Gamma(\Upsilon)$ taking each time the edge with the
positive direction. The final vertex $A$ is joined with some vertex $C$ by
two edges labelled $e$. $\phi$ passes them both in the positive direction,
and then returns to $B$ taking the same edges as for the first time.

\end{enumerate}

The automorphism group $\name{Aut}(F_N)$ acts on the set of subgroups of
$F_N$. The correspondence of Theorem~\ref{Th:SubgrCov} allows to consider
this action as an action on the set of coverings of the union ${\cal U}$,
with a marked vertex chosen in each covering.

\begin{theorem}\label{Th:ActSubgr}
For any tree-like covering $\Upsilon$ and any element $s \in BC_N \subset
\name{Aut}(F_N)$ the covering $s(\Upsilon)$ is tree-like, too. Actions of
the group $BC_N$ on the set of coverings and on the set of trees are
compatible: $\Gamma(s(\Upsilon)) = s(\Gamma(\Upsilon))$.
\end{theorem}

Note that $\Gamma(\Upsilon)$ is a tree with one marked vertex. The group
$BC_N$ acts on such trees, too: the action of Subsection~\ref{SSec:ActTree}
preserves the set of vertices of the tree, so that we can assume that the
marked vertex simply remains unchanged.

\begin{proof}
Let $\Upsilon'$ be a covering such that $\Gamma(\Upsilon') =
s(\Gamma(\Upsilon))$. To prove the statement we must show that the elements
$s(a_{\nu,e}^{\Gamma(\Upsilon)})$ and $s(b_{\nu,e}^{\Gamma(\Upsilon)})$
belong to the subgroup $\Group(\Upsilon')$. Since the subgroups
$s(\Group(\Upsilon))$ and $\Group(\Upsilon')$ have the same index, they
coincide, and thus $s(\Upsilon) = \Upsilon'$.

It is enough to consider the case when $s$ is either $\lambda$ or some
$u_k$. For brevity we will show in detail why
$s(a_{\nu,e}^{\Gamma(\Upsilon)} \in \Group(\Upsilon')$; the proof for
$b_{\nu,e}$ is similar.

{\def \theenumi {\arabic{enumi}}
\def \labelenumi {Case~\theenumi.}
\expandafter\def\csname p@enumi\endcsname {}
\begin{enumerate}
\item\label{It:Lambda} $s = \lambda$. Here $s(a_{\nu,e}^{\Gamma(\Upsilon)})
= a_{\nu,s(e)}^{s(\Gamma(\Upsilon))} \in \Group(\Upsilon')$.

\item\label{It:NAdj} $s = u_k$, and the edges $e_{k-1}$ and $e_k$ in
$\Gamma(\Upsilon)$ are not adjacent. The following situations may occur:

{\def \theenumii {\arabic{enumi}.\arabic{enumii}}
\def \labelenumii {\theenumii.}
\expandafter\def\csname p@enumii\endcsname {}
\begin{enumerate}
\item\label{It:NAdjNotK} $e$ is neither $e_k$ nor $e_{k-1}$. In this case
$s(a_{\nu,e}^{\Gamma(\Upsilon)}) = s(\tau_\nu) e s(\tau_\nu)^{-1}$. The
word $s(\tau_\nu)$ differs from $\tau_\nu$ by substitutions $e_{k-1}
\mapsto e_k,\, e_k \mapsto e_k^{-1} e_{k-1} e_k$. Let $C, D$ be the ends of
the edge $e_{k-1}$, and $E,F$ be the ends of the edge $e_k$ in the graph
$\Gamma(\Upsilon)$. Then in $s(\Gamma(\Upsilon))$ the edge $CD$ is labelled
$e_k$, and $EF$ is labelled $e_{k-1}$. Since the vertices $C, D, E, F$ are
all different, the graph $\Upsilon'$ contains loops $\ell_E$ and $\ell_F$
labelled $e_k$ and attached to vertices $E$ and $F$, respectively.  Thus,
$s(\tau_\nu)$ is just $\tau_{\nu'}$ (in the graph $\Upsilon'$) where $\nu'$
coincides with $\nu$ except for the edge $EF$ which is replaced by the edge
sequence $\ell_E^{-1},\, EF,\, \ell_F$. The last edge of the path $\nu$ in
$\Gamma(\Upsilon)$ is not $e$, and thus, the last edge of $\nu'$ in
$\Gamma(\Upsilon')$ is not $e$ either, and therefore in $\Upsilon'$ there
exists a loop labelled $e$ and attached to the last vertex of $\nu'$. So,
$s(a_{\nu,e}^{\Gamma(\Upsilon)}) \in \Group(\Upsilon')$.

\item\label{NAdjK-1} $e = e_{k-1}$. In this case
$s(a_{\nu,e}^{\Gamma(\Upsilon)}) = s(\tau_\nu) e_k s(\tau_\nu)^{-1}$, and
the subsequent proof is exactly the same as in the
case~\ref{It:NAdjNotK}.

\item\label{It:NAdjKNotK-1} $e = e_k$, and the last edge of the path $\nu$
is not $e_{k-1}$. In this case $s(a_{\nu,e}^{\Gamma(\Upsilon)}) =
s(\tau_\nu) e_k^{-1} e_{k-1} e_k s(\tau_\nu)^{-1}$. The element
$s(\tau_\nu)$ is $\tau_{nu'}$ where $\nu'$ is exactly as in the
case~\ref{It:NAdjNotK}. The last edge of the path $\nu'$ in $\Upsilon'$ is
labelled neither $e_k$ nor $e_{k-1}$, and thus there exist loops $\ell_k$
and $\ell_{k-1}$ with these labels attached to the last vertex of $\nu'$.
So, $s(a_{\nu,e}^{\Gamma(\Upsilon)}) = \tau_{\nu'} \tau_\phi
\tau_{\nu'}^{-1}$ where $\phi$ is $\ell_k^{-1} \ell_{k-1} \ell_k$.
Therefore $s(a_{\nu,e}^{\Gamma(\Upsilon)}) \in \Group(\Upsilon')$.

\item\label{It:NAdjKK-1} $e = e_k$, and the last edge of the path $\nu$ is
$e_{k-1}$. Thus the last edge (say, $EF$) of the path $\nu'$ in $\Upsilon'$
is $e_k$, and there is a loop $\ell$ labelled $e_{k-1}$ and attached to
$E$. The subsequent proof is just the same as in the
case~\ref{It:NAdjKNotK-1} but the path $\phi$ is taken to be $FE,\, \ell,\,
EF$.
\end{enumerate}
}

\item\label{It:Adj} $s = u_k$, and the edges $e_{k-1}$ and $e_k$ in
$\Gamma(\Upsilon)$ are adjacent: $e_{k-1} = CD$, $e_k = DE$. The following
situations may occur here:
{\def \theenumii {\arabic{enumi}.\arabic{enumii}}
\def \labelenumii {\theenumii.}
\expandafter\def\csname p@enumii\endcsname {}
\begin{enumerate}
\item\label{It:NotPass} The path $\nu$ passes neither $e_k$ nor $e_{k-1}$.
In this case $s(\tau_\nu) = \tau_\nu$, and $s(e)$ is either $e$ (if $e \ne
e_k, e_{k-1}$), or $e_k$ (if $e = e_{k-1}$), or $e_k^{-1} e_{k-1} e_k$ (if
$e = e_k$). In all cases the loops labelled $e$, $e_k$ and $e_{k-1}$ are
attached to the final vertex $F$ of the path $\nu$ (both in $\Upsilon$ and
$\Upsilon'$ because the edges of $\nu$ does not change their labels), and
$s(e)$ is $\tau_\phi$ for some path $\phi$ starting and ending at $F$.
Thus $s(a_{\nu,e}^{\Gamma(\Upsilon)}) = \tau_\nu \tau_\phi \tau_\nu^{-1} \in
\Group(\Upsilon')$. The same reasoning applies to the case when $\nu$
passes $e_{k-1}$ (the edge $CD$) but not $e_k$.

\item\label{It:PassK} $\nu$ passes $e_k$ (the edge $DE$) but not $e_{k-1}$.
Graph $\Upsilon'$ does not contain the edge $DE$ but contains the edge $CE$
labelled $e_{k-1}$ (and the edge $CD$ labelled $e_k$). Thus there exists a
loop $\ell$ (in $\Upsilon$) attached to $E$ and labelled $e_k$. So,
$s(\tau_\nu) = \tau_{\nu'}$ where the path $\nu'$ coincides with $\nu$
everywhere except the edge $DE$ which is replaced by the edge sequence
$DC,\, CE,\, \ell$. The word $s(e)$ is handled exactly like in the
case~\ref{It:NotPass}. Thus, $s(a_{\nu,e}^{\Gamma(\Upsilon)}) = \tau_{\nu'}
s(e) \tau_{\nu'}^{-1} \in \Group(\Upsilon')$.

\item\label{It:PassBoth} $\nu$ passes both $e_{k-1}$ and $e_k$. Since
$\Gamma(\Upsilon)$ is a tree, it passes these edges subsequently, either as
$CD,\, DE$, or as $ED,\, DC$. In the first case the reasoning is the same
as in the case~\ref{It:PassK} but $\nu'$ is obtained from $\nu$ by a change
of the segment $CD,\, DE$ to the segment $CE,\, \ell$. The second case is
similar.\qed
\end{enumerate}
}
\end{enumerate}
}
\end{proof}


\begin{thebibliography}{99}

\bibitem{Alexander} Alexander, J.~W. Deformations of $n$-cell. {\em Proc.
Nat. Acad. Sci. USA}, {\bf 9} (1923), pp.~406--407.

\bibitem{AniLan} Anisov, S.~S.; Lando, S.~K. Topological complexity of
$\Torus^2$-bundles over the circle. {\em Amer. Math. Soc. Transl.}, {\bf
185} (1998), No.~2.

\bibitem{Arn} Arnold, V.~I. Critical points of functions and the
classification of caustics, {\em Uspekhi~Mat.~Nauk}, {\bf 29} (1974), no.3,
pp.~243--244 (in Russian).

\bibitem{Birman} Birman, J.~S. {\em Braids, links, and mapping class
groups}. Annals of Mathematics Studies, No. 82. Princeton University Press,
Princeton, N.J.; University of Tokyo Press, Tokyo, 1974. ix+228 pp.

\bibitem{Cohen} Cohen, F.~R. On the mapping class groups for punctured
spheres, the hyperelliptic mapping class groups, $SO(3)$ and $Spin^c(3)$,
{\em Amer.\ J.\ Math.} {\bf 115} (1993), pp.~389--434.

\bibitem{FuFo} Fomenko, A.~T., Fuchs, D.~B. {\em A course in homotopic
topology}, Nauka Publishers, Moscow, 1989. 494 pp.

\bibitem{Thurston} Hatcher, A., Thurston, W. A presentation for mapping
class group of a closed orientable surface, {\em Topology}, {\bf 19}
(1980), No.~3, pp.~221--237.

\bibitem{Looij} Looijenga, E.~J.~N. The complement of the bifurcation
variety of a simple singularity, {\em Invent.~Math.}, {\bf 23} (1974),
pp.~105--116.

\bibitem{Lyash} Lyashko, O.~V. Geometry of bifurcation diagrams, {\em
J.~Soviet Math.}, {\bf 27} (1984), pp.~2736--2759.

\bibitem{PresBraid} Sergiescu, V. Graphes planaires et presentations des
groupes de tresses, {\em Math.~Z.}, {\bf 214} (1993), pp.~477--490.

\end{thebibliography}
\end{document}